\documentstyle[11pt]{article}
\begin{document}
\author{{Shiqiu Zheng$^{1, 2}$\thanks{E-mail: shiqiumath@163.com (S. Zheng).}\ , \ \ Shoumei Li$^{1}$\thanks{Corresponding author, E-mail: li.shoumei@outlook.com (S. Li).}}
  \\
\small(1, College of Applied Sciences, Beijing University of Technology, Beijing 100124, China)\\
\small(2, College of Sciences, North China University of Science and Technology, Tangshan 063009, China)\\
}
\date{}
\title{\textbf{Representation theorem for generators of quadratic BSDEs}\thanks{This work is supported by the National Natural Science Foundation of China (No. 11171010) and Natural Science Foundation of Beijing (No. 1132008).}}\maketitle

\textbf{Abstract:}\quad  In this paper, we establish a general representation theorem for generator of backward
stochastic differential equation (BSDE), whose generator has a quadratic growth in $z$. As some applications, we obtain a general converse comparison theorem of such quadratic BSDEs and uniqueness theorem, translation invariance for quadratic $g$-expectation.
\\



\section{Introduction}
It is well-known that Pardoux and Peng (1990) firstly establishes the
existence and uniqueness of solution of the nonlinear BSDE under
the Lipschitz condition. Since then, the theory of BSDEs
has gone through rapid development in many different
areas such as probability, PDE, stochastic control, and mathematical finance, etc. An important study in BSDEs theory is to interpret the relation between the generators and solutions of BSDEs. In this topic, the representation theorem for generator of BSDE is established, which shows a relation between
generators and solutions of BSDEs in limit form. It plays an important role in the study in BSDEs theory and nonlinear expectation theory.

Representation theorem of generator is firstly established by Briand et al. (2000) for BSDEs with Lipschitz generators under two additional assumptions that $E[\sup_{0\leq t\leq T}|g(t,0,0)|^2] <\infty$ and $(g(t, y, z))_{t\in[0,T]}$ is continuous
in $t$, in order to study the converse comparison theorem of BSDEs. After a series of studies by Jiang (2005a, 2005b), two assumptions mentioned above in Briand et al. (2000) are eliminated by Jiang (2005c, 2008). Since then, the representation theorem for generators is further studied for BSDEs with linear growth (\emph{\textbf{generator has a linear growth in $z$}}). One can see Fan and Jiang (2010), Fan et al. (2011) and Jia (2008), etc.

Quadratic BSDEs (\emph{\textbf{generator has a quadratic growth in $z$}}) have been firstly studied by Kobylanski (2000), then by many papers. Such BSDEs have many important applications in PDE, stochastic control and finance. To interpret the relation between the generators and solutions of quadratic BSDEs considered in Kobylanski (2000). Ma and Yao (2010, Theorem 4.1) establish a representation theorem for generator of BSDE when generator satisfies some conditions similar as in Kobylanski (2000) and two additional assumptions (see conditions ($g1$) and ($g2$) in Ma and Yao (2010, Theorem 4.1)). Recently, Jia and Zhang (2015, Theorem 3.1) establish a representation theorem for generator of BSDE with so called Lipschitz-quadratic generator  (see assumption (H) in Jia and Zhang (2015)). The main result of this paper is a general representation theorem for generator of BSDE, whose generator only satisfies the conditions in Kobylanski (2000). We firstly obtain this representation theorem in almost surely sense, then in $L^p$ sense for $p>0$. This phenomenon is different from the representation theorems for BSDEs with linear growth, which are obtained in $L^p$ or ${\cal{H}}^p$ sense for $1\leq p<2$ or $p=2$ (see Fan and Jiang (2010), Fan et al. (2011), Jia (2008) and Jiang (2008), etc). Recently, some results in this paper has been used to establish a representation theorem for BSDE whose generator is monotonic and convex growth in $y$ and quadratic in $z$ (see Zheng and Li (2015)). A discussion on the differences between the present
work and some known results is given in Remark 3.7.

The famous comparison theorem of BSDEs shows that we can compare the solutions of BSDEs through comparing generators. Converse comparison theorem for BSDEs shows that we can compare the generators through comparing the solutions of BSDEs. It is firstly studied by Chen (1997), then by Briand et al. (2001), Coquet et al.(2001), Fan et al. (2011) and Jiang (2004, 2005b, 2005c), etc., for BSDEs with linear growth, and by Ma and Yao (2010) for quadratic BSDEs. Using representation theorem obtained in this paper, we establish a general converse comparison theorem for quadratic BSDEs.

The notion of $g$-expectation is introduced by Peng (1997), which is a nonlinear expectation induced by BSDE with Lipschitz generator.
Ma and Yao (2010) consider $g$-expectation induced by quadratic BSDE, called quadratic $g$-expectation, and study its properties. Quadratic $g$-expectation is also studied in Hu et al. (2008). Using representation theorems obtained in this paper, we study uniqueness theorem, translation invariance for quadratic $g$-expectation.

This paper is organized as follows. In section 2, we will present some basic assumptions. In section 3, we will prove some important lemmas and establish representation theorems for quadratic BSDE. In section 4 and section 5, we will give some applications of representation theorem for quadratic BSDEs and quadratic $g$-expectation, respectively.

\section{Preliminaries}
Let $(\Omega ,\cal{F},\mathit{P})$ be a complete probability space
carrying a $d$-dimensional standard Brownian motion ${{(B_t)}_{t\geq
0}}$, starting from $B_0=0$, let $({\cal{F}}_t)_{t\geq 0}$ denote
the natural filtration generated by ${{(B_t)}_{t\geq 0}}$, augmented
by the $\mathit{P}$-null sets of ${\cal{F}}$, let $|z|$ denote its
Euclidean norm, for $\mathit{z}\in {\mathbf{R}}^n$, let $T>0$ be a given real number. We define the following usual
spaces:

${\cal{L}}^p[0,T]=\{f(t):$ Lebesgue measurable function; $\int_0^T|f(t)|^pdt<\infty\},\ \ \ p\geq1; $

$L^p({\mathcal {F}}_T)=\{\xi:\ {\cal{F}}_T$-measurable
random variable; $\|\xi\|_{L^p}=\left({\mathbf{E}}\left[|\xi|^p\right]\right)^{\{\frac{1}{p}\wedge1\}}<\infty\},\ \ \ p>0; $

$L^\infty({\mathcal {F}}_T)=\{\xi:{\cal{F}}_T$-measurable
random variable; $\|\xi\|_\infty=\textrm{esssup}_{\omega\in\Omega}|\xi|<\infty\};$

${\cal{H}}^p_T({\mathbf{R}}^d)=\{\psi:$ predictable
process;
$\|\psi\|_{{\cal{H}}^p}=\left({\mathbf{E}}\left[\int_0^T|\psi_t|^pdt\right]\right)^{\{\frac{1}{p}\wedge1\}}
<\infty \},\ \ \ \ p>0;$

${\cal{H}}^\infty_T({\mathbf{R}}^d)=\{\psi:$ predictable
process;
$\|\psi\|_{{\cal{H}}^\infty}=\textrm{esssup}_{(w,t)\in\Omega\times[0,T]}|\psi_t|<\infty \};$

${\mathcal{S}}^\infty
 _T({\mathbf{R}})=\{\psi:$ continuous predictable
process in ${\cal{H}}^\infty_T({\mathbf{R}})\}.$

Let us consider a function $g$
$${g}\left( \omega ,t,y,z\right) : \Omega \times [0,T]\times \mathbf{%
R\times R}^{\mathit{d}}\longmapsto \mathbf{R}$$ such that
$\left(g(t,y,z)\right)_{t\in [0,T]}$ is progressively measurable for
each $(y,z)\in\mathbf{%
R\times R}^{\mathit{d}}$ and $g$ always takes the following form:
$$g(t,y,z)=g_1(t,y,z)y+g_2(t,y,z).$$
we give the following assumptions for $g$:

\begin{itemize}
\item (A1). For each $(y,z)\in{\mathbf{R}}\times{\mathbf{R}}^d,$ both $g_1$ and $g_2$ are progressively measurable and for each $(t,w)\in [0,T]\times\Omega,$ both $g_1$ and $g_2$ are continuous in $(y,z);$
\item (A2). There exist constants $\alpha,\beta,b$ and a continuous increasing function $l:\ \textbf{R}^+\rightarrow\textbf{R}^+,$ such that for $(t,y,z)\in[0,T]\times\mathbf{%
R\times R}^{\mathit{d}},$
$$P-a.s.,\ \ \beta\leq g_1(t,y,z)\leq \alpha\ \ \ \ \textrm{and}\ \ \ \ \ |g_2(t,y,z)|\leq b+l(|y|)|z|^2;$$
 \item (A3). For any $M>0,$ there exist a function $k(t)\in {\cal{L}}^2[0,T]$ and a constant $C$ such that
$\forall (t,y,z)\in[0,T]\times[-M,M]\times\mathbf {R}^{\mathit{d}},$
$$P-a.s.,\ \ \left|\frac{\partial g}{\partial z}(t,y,z)\right|\leq k(t)+C|z|;$$
  \item (A4). For any $\varepsilon>0,$ there exists a function $h_{\varepsilon}(t)\in {\cal{L}}^1[0,T],$ such that for $(t,y,z)\in[0,T]\times\mathbf{%
R\times R}^{\mathit{d}},$
$$P-a.s.,\ \ \frac{\partial g}{\partial y}(t,y,z)\leq h_{\varepsilon}(t)+\varepsilon|z|^2.$$
\end{itemize}
\textbf{Remark 2.1} If $g$ satisfies assumption (A2), then it also satisfies the following (A2)$^\ast$.
\begin{itemize}
  \item (A2)$^\ast$. For any constant $M>0,$ there exists a constant
$\lambda_M=\max\{|\alpha|, |\beta|, |b|,l(M)\}$  such that
$\forall (t,y,z)\in[0,T]\times[-M,M]\times\mathbf {R}^{\mathit{d}},$
$$P-a.s.,\ \ |{g}(t,y,z)|\leq \lambda_M(1+|y|+|z|^2).$$
\end{itemize}

In this paper, we consider the following BSDE introduced by Pardoux and Peng (1990):
$$Y_t=\xi +\int_t^Tg(s,Y_s,Z_s)
ds-\int_t^TZ_sdB_s,\ \ \ t\in[0,T]\eqno(1)$$
where $g$ is called generator, $\xi$ and $T$ are called terminal variable and terminal time, respectively. The BSDE (1) is usually called BSDE with parameter $(g,T,\xi).$

Now, we introduce a stochastic differential equation (SDE). Suppose $b(\cdot,\cdot,\cdot):\Omega\times [0,T]\times \textbf{R}^m\mapsto \textbf{R}^m$ and $\sigma(\cdot,\cdot,\cdot):\Omega\times [0,T]\times \textbf{R}^m\mapsto \textbf{R}^{m\times d}$ satisfy the following two conditions:
\begin{itemize}
  \item (H1). There exists a constant $\mu\geq0$ such that $P$-$a.s.$, $ \forall t\in[0,T],\
\forall x,y\in \textbf{R}^{\mathit{m}},$
$$|{b}(t,x)-b(t,y) |+|\sigma(t,x)-\sigma(t,y) |\leq \mu|x-y|.$$
  \item (H2). There exists a constant $\nu\geq0$ such that $P$-$a.s.$, $ \forall t\in[0,T],\
\forall x\in \textbf{R}^{\mathit{m}},$
$$|b(t,x)|+|\sigma(t,x)|\leq \nu\left(1+|x|\right).$$
\end{itemize}
Given $(t,x)\in[0,T[\times \textbf{R}^m$, by SDE theory, the following SDE:
$$\left\{
  \begin{array}{ll}
     X_s^{t,x}=x +\int_t^sb(u,X_u^{t,x})du+\int_t^s\sigma(u,X_u^{t,x})dB_u,\ \ \  s\in]t,T],\\
     X_s^{t,x}=x\ \ \   s\in[0,t],
  \end{array}
\right.$$
have a unique continuous adapted solution $X_s^{t,x}$.

In the end of this section, we introduce the following Lebesgue Lemma, which plays an important role in this paper.\\\\
\
\textbf{Lemma 2.2} (Hewitt and Stromberg (1978, Lemma 18.4)) Let $f\in {\cal{L}}^1[0,T]$. Then for almost every $t\in[0,T[$, we have
$$\lim\limits_{\varepsilon\rightarrow0^+}\frac{1}{\varepsilon}
\int_{t}^{t+\varepsilon}|f(u)-f(t)|ds=0.$$

\section{Representation theorems for quadratic BSDEs}
In this section, we will study the representation theorem for generators of quadratic BSDEs. Firstly, we recall some known results in Kobylanski (2000). By Kobylanski (2000, Theorem 2.3), if $g$ assumptions (A1) and (A2), for stopping time $\delta\leq T,\ a.s.$ and $\xi\in L^\infty({\mathcal {F}}_\delta),$ BSDE with parameter $(g,\delta,\xi)$ has a maximal solution $(\overline{{Y}}_t,\overline{{Z}}_t)\in{\mathcal{S}}^\infty_\delta
({\mathbf{R}})\times{\cal{H}}^{2}_\delta({\mathbf{R}}^d)$ and a minimal solution $(\underline{{Y}}_t,\underline{{Z}}_t)\in{\mathcal{S}}^\infty_\delta
({\mathbf{R}})\times{\cal{H}}^{2}_\delta({\mathbf{R}}^d)$ in the sense that for any solution $({Y}_t,{Z}_t)$ of BSDE with parameter $(g,\delta,\xi)$, we have $\underline{{Y}}_{t\wedge\delta}\leq{Y}_{t\wedge\delta}\leq\overline{{Y}}_{t\wedge\delta},\ a.s.,$ for each $t\in[0,T].$ Moreover, we have the following two facts:
$$\|Y_{t\wedge\delta}\|_{{\cal{H}}^\infty}\leq(\|\xi\|_\infty+|b|T)e^{\alpha^+ T},\eqno(2)$$
and for each $\varepsilon\in]0,T-t]$ and stopping time $\tau\in ]0,T-t],$ if $({Y}_{s}^{t+\varepsilon\wedge\tau},{Z}_{s}^{t+\varepsilon\wedge\tau})$ is a solution of BSDE with parameter $(g,t+\varepsilon\wedge\tau,0),$ then
$$P-a.s.,\ \ \sup_{t\leq s\leq t+\varepsilon\wedge\tau}|Y_{s}^{t+\varepsilon\wedge\tau}|\leq |b|\varepsilon e^{\alpha^+\varepsilon}.\eqno(3)$$
In fact, since $\|{Y}_{t\wedge\delta}\|_{{\cal{H}}^\infty}<\infty$, we define a continuous function $\phi({y}):\textbf{R}\rightarrow[0,1],$ such that
$$\phi({y})
:=\left\{
 \begin{array}{ll}
  1, \ \ \ \ |{y}|\leq \|{Y}_{t\wedge\delta}\|_{{\cal{H}}^\infty},\\
  0,\ \ \ \ |{y}|> \|{Y}_{t\wedge\delta}\|_{{\cal{H}}^\infty}+1.
  \end{array}
  \right.
$$
Set
$$\hat{g}(r,y,z):={g}_1(r,y,z)y+\phi(y){g}_2(r,y,z),$$
then by (A2), we can check that $\hat{g}(r,y,z)$ satisfies condition (H0) in Kobylanski (2000, Proposition 2.1 and Corollary 2.2). Clearly, $({Y}_t,{Z}_t)$ is also a solution of BSDE with parameter $(\hat{g},\delta,\xi)$. Then we can get (2) from Kobylanski (2000, Corollary 2.2), immediately. Similarly, we also can get (3) from Kobylanski (2000, Proposition 2.1), immediately.

Now, we will prove three lemmas in the following.\\\\
\
\textbf{Lemma 3.1} Let $g$ satisfy (A1) and (A2). Then for each $t\in[0,T[$ and stopping time $\tau\in ]0,T-t],$ there exists a constant $\gamma>0$ such that for each $\varepsilon\leq \gamma,$
$$E\left[\int_{t}^{t+\varepsilon\wedge\tau}|{Z}_r^{t+\varepsilon\wedge\tau}|^2dr|{\cal{F}}_t\right]\leq C\varepsilon^2,$$
where $(Y_{s}^{t+\varepsilon\wedge\tau},Z_{s}^{t+\varepsilon\wedge\tau})$ is an arbitrary solution of BSDE with parameter $(g,t+\varepsilon\wedge\tau,0)$ and $C$ is a constant only depending on $\alpha, \beta, b,T,\gamma.$\\\\
\
\emph{Proof.} For $t\in[0,T[,\ \varepsilon\in]0,T-t]$ and stopping time $\tau\in]0,T-t],$ we consider the following BSDEs with parameter $(g,t+\varepsilon\wedge\tau,0)$
$${Y}_s^{t+\varepsilon\wedge\tau}=\int_{s}^{{t+\varepsilon\wedge\tau}}g(r,{Y}_r^{t+\varepsilon\wedge\tau},
{Z}_r^{t+\varepsilon\wedge\tau})dr-\int_{s}^{t+\varepsilon\wedge\tau}{Z}_r^{t+\varepsilon\wedge\tau}dB_r.$$
By (2), we have
$$P-a.s., \ \ \sup_{0\leq s\leq t+\varepsilon\wedge\tau}|{Y}_{s}^{t+\varepsilon\wedge\tau}|\leq |b|Te^{\alpha^+T}.$$
We set $\hat{M}:=|b|Te^{\alpha^+T}.$
Applying It\^{o} formula to $|{Y}_s^{t+\varepsilon\wedge\tau}|^2$ for $s\in[t,t+\varepsilon\wedge\tau],$ and in view of (A2)$^\ast$, we can get there exists a constant $\lambda_{\hat{M}}=\max\{|\alpha|, |\beta|, |b|, l(\hat{M})\}$ such that
\begin{eqnarray*}
&&|{Y}_t^{t+\varepsilon\wedge\tau}|^2+\int_{t}^{t+\varepsilon\wedge\tau}|{Z}_r^{t+\varepsilon\wedge\tau}|^2dr\\
&\leq&2\int_{t}^{t+\varepsilon\wedge\tau}|{Y}_r^{t+\varepsilon\wedge\tau}||g(r,{Y}_r^{t+\varepsilon\wedge\tau},
{Z}_r^{t+\varepsilon\wedge\tau})|dr-2\int_{t}^{t+\varepsilon\wedge\tau}{Y}_r^{t+\varepsilon\wedge\tau}{Z}_r^{t+\varepsilon\wedge\tau}dB_r\\
&=&2\lambda_{\hat{M}}\int_{t}^{t+\varepsilon\wedge\tau}|{Y}_r^{t+\varepsilon\wedge\tau}|dr
+2\lambda_{\hat{M}}\int_{t}^{t+\varepsilon\wedge\tau}|{Y}_r^{t+\varepsilon\wedge\tau}|^2dr
+2\lambda_{\hat{M}}\int_{t}^{t+\varepsilon\wedge\tau}|{Y}_r^{t+\varepsilon\wedge\tau}||{Z}_r^{t+\varepsilon\wedge\tau}|^2dr\\
&&-2\int_{t}^{t+\varepsilon\wedge\tau}{Y}_r^{t+\varepsilon\wedge\tau}{Z}_r^{t+\varepsilon\wedge\tau}dB_r\ \ \ \ \ \ \ \ \ \ \ \ \ \ \ \ \ \ \ \ \ \ \ \ \ \ \ \ \ \ \ \ \ \ \ \ \ \ \ \  \ \ \ \ \ \ \ \ \ \ \ \ \ \ \  \ \ \ \ \ \ \ \ \ \ \ \ \  \ \ \ \ \ \ \  (4)
\end{eqnarray*}
Then by (3), we can select $\gamma$ small enough such that for each $\varepsilon\leq \gamma,$ we have
$$\sup_{t\leq s\leq t+\varepsilon\wedge\tau}|Y_s^{t+\varepsilon\wedge\tau}|\leq \frac{1}{4\lambda_{\hat{M}}}.$$
Then by this, (3) and (4), for each $\varepsilon\leq \gamma,$ we have
\begin{eqnarray*}
&&E\left[\int_{t}^{t+\varepsilon\wedge\tau}|{Z}_r^{t+\varepsilon\wedge\tau}|^2dr|{\cal{F}}_t\right]\\
&\leq& 2\lambda_{\hat{M}} E\left[\int_{t}^{t+\varepsilon\wedge\tau}|{Y}_r^{t+\varepsilon\wedge\tau}|dr
+\int_{t}^{t+\varepsilon\wedge\tau}|{Y}_r^{t+\varepsilon\wedge\tau}|^2dr
+\int_{t}^{t+\varepsilon\wedge\tau}|{Y}_r^{t+\varepsilon\wedge\tau}||{Z}_r^{t+\varepsilon\wedge\tau}|^2dr|{\cal{F}}_t\right]\\
&\leq&2\varepsilon^2\lambda_{\hat{M}}^2e^{\varepsilon\lambda_{\hat{M}}}
+2\varepsilon^3\lambda_{\hat{M}}^3e^{2\varepsilon\lambda_{\hat{M}}}
+\frac{1}{2}E\left[\int_{t}^{t+\varepsilon\wedge\tau}|{Z}_r^{t+\varepsilon\wedge\tau}|^2dr|{\cal{F}}_t\right].
\end{eqnarray*}
From above inequality, the proof is complete. $\Box$\\

Inspired by Fan and Jiang (2010, Proposition 3) and Fan et al. (2011, Lemma 3 and Proposition 2), we have the following Lemma 3.2.\\\\
\
\textbf{Lemma 3.2} Let $g$ satisfy (A1) and (A2). For any constant $M>0$ and $(y,x,q)\in[-M,M]\times {\mathbf{R}}^{\mathit{m}}\times {\mathbf{R}}^{\mathit{m}}$, there exists a non-negative process sequence $\{(\psi^n_t)_{t\in[0,T]}\}_{n=1}^{\infty}$ depending on $(y,x,q),$ which satisfies
$$\|\psi^n_t\|_{{\cal{H}}^\infty}\leq C,\ \ \textmd{and}\ \ \forall t\in[0,T],\  \lim\limits_{n\rightarrow\infty}\psi^n_t=0,\ \ P-a.s.,\eqno(5)$$
such that for each $n\in \textbf{N}$ and $(t,\bar{y},\bar{z},\bar{x})\in [0,T]\times[-M,M]\times {\mathbf{R}}^{d+m}$, we have
$$|g(t,\bar{y},\bar{z}+\sigma^\ast(t,\bar{x})q)-g(t,y,\sigma^\ast(t,x)q)|\leq 2n\tilde{\lambda}_M(|\bar{y}-y|+|\bar{z}|^2+|\bar{x}-x|^2)+{\psi}^n_t.$$
In the above, $\tilde{\lambda}_M=4\lambda_M (1+|q|^2|\nu|^2),$ $\lambda_M=\max\{|\alpha|, |\beta|, |b|,l(M)\}$ and $C$ is a constant only dependent on $\tilde{\lambda}_M$ and $(y,x).$\\\\
\
\emph{Proof.} For any constant $M>0$ and $(y,x,q)\in[-M,M]\times {\mathbf{R}}^{\mathit{m}}\times {\mathbf{R}}^{\mathit{m}}$, we set $$f(t,y,z,x):=g(t,y,z+\sigma^\ast(t,x)q).$$
Then by (A2)$^\ast$ and (H2), there exists a constant
$\lambda_M=\max\{|\alpha|, |\beta|, |b|,l(M)\}$ such that
\begin{eqnarray*}
|f(t,y,z,x)|&\leq& \lambda_M (1+|y|+|z+\sigma^\ast(t,x)q|^2)\\
&\leq&\lambda_M (1+|y|+2|z|^2+4|q|^2|\nu|^2(1+|x|^2))\\
&\leq&4\lambda_M (1+|q|^2|\nu|^2)(1+|y|+|z|^2+|x|^2).
\end{eqnarray*}
we set $\tilde{\lambda}_M:=4\lambda_M (1+|q|^2|\nu|^2)$ and
$${f}_n^1(t,y,z,x):=\sup_{(u,v,w)\in \{\textbf{Q}\cap[-M,M]\}\times\textbf{Q}^{d+n}}\{f(t,u,v,w)-2n\tilde{\lambda}_M(|u-y|+|v-z|^2+|w-x|^2)\}.$$
$${f}_n^2(t,y,z,x):=\inf_{(u,v,w)\in \{\textbf{Q}\cap[-M,M]\}\times\textbf{Q}^{d+n}}\{f(t,u,v,w)+2n\tilde{\lambda}_M(|u-y|+|v-z|^2+|w-x|^2)\}.$$
Then by argument of Fan et al. (2011, Lemma 3), we can deduce the following fact:

(i)$|{f}_n^i(t,y,z,x)|\leq2\tilde{\lambda}_M(1+|y|+|z|^2+|x|^2),\ \ \ i=1,2;$

(ii)${f}_n^1(t,y,z,x)\searrow$ and ${f}_n^2(t,y,z,x)\nearrow,$ as $n\longrightarrow\infty;$

(iii)${f}_n^i(t,y,z,x)\longrightarrow f(t,y,z,x)$ as $n\longrightarrow\infty,\ \ \ i=1,2.$\\
Setting
$${\psi}^1_n(t):={f}_n^1(t,y,0,x)\ \  \textmd{and} \ \ {\psi}^2_n(t):={f}_n^2(t,y,0,x).$$
By (i), we have
$$|{\psi}^1_n(t)|\leq2\tilde{\lambda}_M(1+|y|+|x|^2)\ \ \textmd{and}\ \ |{\psi}^2_n(t)|\leq2\tilde{\lambda}_M(1+|y|+|x|^2).$$
Then by (iii), we have
$$\lim_{n\rightarrow\infty}{\psi}^1_n(t)=\lim_{n\rightarrow\infty}{\psi}^2_n(t)=f(t,y,0,x).$$
By the definition of $f^1_n$ and $f^2_n$, we also have
$$f(t,\bar{y},\bar{z},\bar{x})-f(t,y,0,x)\leq 2n\tilde{\lambda}_M(|\bar{y}-y|+|\bar{z}|^2+|\bar{x}-x|^2)+{\psi}^1_n(t)-f(t,y,0,x);$$
$$f(t,\bar{y},\bar{z},\bar{x})-f(t,y,0,x)\geq -2n\tilde{\lambda}_M(|\bar{y}-y|+|\bar{z}|^2+|\bar{x}-x|^2)+{\psi}^2_n(t)-f(t,y,0,x);$$
By setting
$${\psi}^n_t:=|{\psi}^1_n(t)-f(t,y,0,x)|+|{\psi}^2_n(t)-f(t,y,0,x)|,$$
we have
$$|f(t,\bar{y},\bar{z},\bar{x})-f(t,y,0,x)|\leq 2n\tilde{\lambda}_M(|\bar{y}-y|+|\bar{z}|^2+|\bar{x}-x|^2)+{\psi}^n_t.$$
The proof is complete.  $\Box$\\\\
\
\textbf{Lemma 3.3} Let $(\varphi_t)_{t\in[0,T]}$ be a $\mathbf{R}$-valued
progressively measurable process such that
$\|\varphi_t\|_{{\cal{H}}^\infty}<\infty,$ then for almost every $t\in[0,T[$ and any stopping time $\tau_t\in]0,T-t],$ we have
$$P-a.s.,\ \ \ \varphi_t=\lim\limits_{\varepsilon\rightarrow0^+}E\left[\frac{1}{\varepsilon}
\int_{t}^{t+\varepsilon\wedge\tau_t}\varphi_rdr|{\cal{F}}_t\right].$$
\emph{Proof.} For $t\in[0,T[,$ let $\tau_t\in]0,T-t]$ is a stopping time. Set $l^\varepsilon_t:=\frac{1}{\varepsilon}\int_{t}^{t+\varepsilon\wedge\tau_t}\varphi_rdr$ for $t\in[0,T[.$ we can easily check $l^\varepsilon_t$ is a continuous process on $[0,T[.$ Thus $l^\varepsilon_t$ is measurable w.r.t. ${\cal{B}}([0,T[)\otimes{\cal{F}}_T.$ we set
$$\Xi_1:=\{(t,\omega)\in[0,T[\times\Omega:\liminf_{\varepsilon\rightarrow0^+} l^\varepsilon_t\neq\varphi_t\},\ \ \textmd{and}\ \ \Xi_2:=\{(t,\omega)\in[0,T[\times\Omega:\limsup_{\varepsilon\rightarrow0^+} l^\varepsilon_t\neq\varphi_t\}.$$
Thus $\Xi_i$ is a measurable set w.r.t. ${\cal{B}}([0,T[)\otimes{\cal{F}}_T.$ For $\omega\in\Omega$ and $t\in[0,T[,$ we set
$$\Xi_i^\omega=\{t:(t,\omega)\in\Xi_i\}, \ \ i=1,2,\ \ \textmd{and}\ \ (\Xi_1\cup\Xi_2)^\omega=\{t:(t,\omega)\in\Xi_1\cup\Xi_2\}.$$
By Lemma 2.2, we have $P-a.s.,$  for almost every $t\in[0,T[,$
$\lim\limits_{\varepsilon\rightarrow0^+} l^\varepsilon_t
=\varphi_t.$ Thus, we have  $P-a.s.,$
$$\Lambda(\Xi_i^\omega)\leq\Lambda((\Xi_1\cup\Xi_2)^\omega)=0,$$
where $\Lambda$ is the Lebesgue measure.
By Fubini theorem, we have $$\Lambda\otimes P(\Xi_i)=\int_\Omega\Lambda(\Xi_i^\omega)dP=0,\ \  i=1,2.$$
Thus $\Lambda\otimes P((\Xi_1\cup\Xi_2))=0.$ Set $\Xi=\Xi_1^c\cap\Xi_2^c,$ we have $\forall(t,\omega)\in\Xi,$ $\lim\limits_{\varepsilon\rightarrow0^+} l^\varepsilon_t
=\varphi_t.$ Furthermore,
$$\Lambda\otimes P(\Xi)=\Lambda\otimes P(\Xi_1^c\cap\Xi_2^c)=\Lambda\otimes P((\Xi_1\cup\Xi_2)^c)=T-\Lambda\otimes P((\Xi_1\cup\Xi_2))=T.$$
Thus for almost every $t\in[0,T[,$ $P-a.s.,$
$\lim\limits_{\varepsilon\rightarrow0^+} l^\varepsilon_t
=\varphi_t.$ Then by Jensen inequality and dominated convergence theorem, we have for almost every $t\in[0,T[,$
\begin{eqnarray*}
E\left|E\left[\lim\limits_{\varepsilon\rightarrow0^+}
 l^\varepsilon_t|{\cal{F}}_t\right]
-\varphi_t\right|
&=&E\left|E\left[\lim\limits_{\varepsilon\rightarrow0^+}
 l^\varepsilon_t-\varphi_t|{\cal{F}}_t\right]\right|\\
&\leq&E\left|\lim\limits_{\varepsilon\rightarrow0^+}
 l^\varepsilon_t
-\varphi_t\right|\\
&=&0.
\end{eqnarray*}
By this and dominated convergence theorem, for almost every $t\in[0,T[$, we have,
$$P-a.s.,\ \ \lim\limits_{\varepsilon\rightarrow0^+}E\left[ l^\varepsilon_t|{\cal{F}}_t\right]
=E\left[\lim\limits_{\varepsilon\rightarrow0^+} l^\varepsilon_t|{\cal{F}}_t\right]=\varphi_t.$$
The proof is complete.  $\Box$\\

The following Theorem 3.4 is a general representation theorem for generators of quadratic BSDEs, which is the main result of this paper.\\\\
\
\textbf{Theorem 3.4} Let $g$ satisfy (A1) and (A2). Then for each $(y,x,q)\in{\mathbf{R}}\times{\mathbf{R}}^{\mathit{m}}\times {\mathbf{R}}^{\mathit{m}}$ and almost every $t\in[0,T[$, we have
$$P-a.s.,\ \ \ g\left(t,y,\sigma^\ast(t,x)q\right)+q\cdot b(t,x)
=\lim\limits_{\varepsilon\rightarrow0^+}\frac{1}{\varepsilon}
\left(Y_t^{t+\varepsilon\wedge\tau}-y\right),$$
where $\tau=\inf\{s\geq0:|X_{t+s}^{t,x}|> C_0\}\wedge (T-t)$ for a constant $C_0>|x|$ and $(Y_s^{t+\varepsilon\wedge\tau},Z_s^{t+\varepsilon\wedge\tau})$ is an arbitrary solution of BSDE with parameter $(g,t+\varepsilon\wedge\tau,y+q\cdot(X_{t+\varepsilon\wedge\tau}^{t,x}-x)).$\\\\
\
\textit{Proof.} Given $(y,x,q)\in{\mathbf{R}}\times{\mathbf{R}}^{\mathit{m}}\times {\mathbf{R}}^{\mathit{m}}$ and a constant $C_0>|x|.$ For $t\in[0,T[,$ we define the following stopping time:
$$\tau:=\inf\left\{s\geq0:|X_{t+s}^{t,x}|> C_0\right\}\wedge (T-t)$$
By the continuity of $X_{t+s}^{t,x}$, we have
$0<\tau\leq T-t$ and $\forall s\in[0,T],$
$$\|X_{t+{s\wedge\tau}}^{t,x}\|_\infty\leq C_0.\eqno(6)$$
For $\varepsilon\in]0,T-t]$, let $\left(Y_s^{t+\varepsilon\wedge\tau},Z_s^{t+\varepsilon\wedge\tau}\right)$ be a solution of BSDE with parameter $(g,t+\varepsilon\wedge\tau,y+q\cdot(X_{t+\varepsilon\wedge\tau}^{t,x}-x))$ and for $s\in[t,{t+\varepsilon\wedge\tau}],$ we set
$$
\tilde{Y}_s^{{t+\varepsilon\wedge\tau}}:=Y_s^{{t+\varepsilon\wedge\tau}}-(y+q\cdot(X_{s}^{t,x}-x)),\ \ \ \
\tilde{Z}_s^{t+\varepsilon\wedge\tau}:=Z_s^{t+\varepsilon\wedge\tau}-\sigma^\ast(t,X_{s}^{t,x})q.\eqno(7)
$$
By (2),  (6) and (7), there exists a constant $\tilde{M}$ depending $b,\alpha,y,x,q,T$ and $C_0,$ such that for $s\in[t,{t+\varepsilon\wedge\tau}],$ we have
$$\|\tilde{Y}_s^{{t+\varepsilon\wedge\tau}}\|_\infty\leq\|{Y}_s^{{t+\varepsilon\wedge\tau}}\|_\infty
 +\|y+q\cdot(X_{s}^{t,x}-x)\|_\infty\leq\tilde{M}.\eqno(8)$$
Applying It\^{o} formula to $\tilde{Y}_s^{t+\varepsilon\wedge\tau}$ for $s\in[t,{t+\varepsilon\wedge\tau}],$ we have
\begin{eqnarray*}
\tilde{Y}_s^{t+\varepsilon\wedge\tau}&=&\int_{s}^{{t+\varepsilon\wedge\tau}}\left(g(r,\tilde{Y}_r^{t+\varepsilon\wedge\tau}+y
+q\cdot(X_{r}^{t,x}-x),
\tilde{Z}_r^{t+\varepsilon\wedge\tau}+\sigma^\ast(r,X_{r}^{t,x})q)+q\cdot b(r,X_{r}^{t,x})\right)dr\\
&&-\int_{s}^{t+\varepsilon\wedge\tau}\tilde{Z}_r^{t+\varepsilon\wedge\tau}dB_r.\ \ \ \ \ \ \ \ \ \ \ \ \ \ \ \ \ \ \ \ \ \ \ \ \ \ \ \ \ \ \ \ \ \ \ \ \ \ \ \ \ \ \ \ \ \ \ \ \ \ \ \ \ \ \ \ \ \ \ \ \ \ \ \ \ \ \ \ \ \ \ \ \ \ \ \ (9)
\end{eqnarray*}
Now, we define a continuous function $\tilde{\phi}(\tilde{y}):\textbf{R}\rightarrow[0,1],$ such that
$$\tilde{\phi}(\tilde{y})
:=\left\{
 \begin{array}{ll}
  1, \ \ \ \ |\tilde{y}|\leq \tilde{M},\\
  0,\ \ \ \ |\tilde{y}|> \tilde{M}+1.
  \end{array}
  \right.
$$
Set
$$\tilde{g}(r,\tilde{y},\tilde{z}):=\tilde{g_1}(r,\tilde{y},\tilde{z})y+\tilde{g_2}(r,\tilde{y},\tilde{z}),$$
where
\begin{eqnarray*}
\tilde{g_1}(r,\tilde{y},\tilde{z})&:=&g_1(r,\tilde{y}+y+q\cdot(X_{r}^{t,x}-x),\tilde{z}+\sigma^\ast(r,X_{r}^{t,x})q),\\
\tilde{g_2}(r,\tilde{y},\tilde{z})
&:=&\tilde{\phi}(\tilde{y})g_2(r,\tilde{y}+y+q\cdot(X_{r}^{t,x}-x),\tilde{z}+\sigma^\ast(r,X_{r}^{t,x})q)+q\cdot b(r,X_{r}^{t,x}).
\end{eqnarray*}
Clearly, $\tilde{g}$ satisfy the assumptions (A1). Moreover, by (A2), (H2), (6) and (8), we have for $r\in[0,t+\varepsilon\wedge\tau],$
\begin{eqnarray*}
\beta\leq\tilde{g}_1(r,\tilde{y},\tilde{z})|\leq \alpha,
\end{eqnarray*}
and
\begin{eqnarray*}
|\tilde{g}_2(r,\tilde{y},\tilde{z})|&=&|\phi(\tilde{y})g_2(r,\tilde{y}+y+q\cdot(X_{r}^{t,x}-x),\tilde{z}+\sigma^\ast(r,X_{r}^{t,x})q)+q\cdot b(r,X_{r}^{t,x})|\\
&&\leq |b|+|\tilde{\phi}(\tilde{y})|l(|\tilde{y}+y+q\cdot(X_{r}^{t,x}-x)|)|\tilde{z}+\sigma^\ast(r,X_{r}^{t,x})q|^2+|q\cdot b(r,X_{r}^{t,x})|
\\
&&\leq \tilde{b}+2l(|\tilde{y}|+\tilde{M})|\tilde{z}|^2,
\end{eqnarray*}
where $\tilde{b}=|b|+2q^2\nu^2(1+C_0)^2l(2\tilde{M}+1)+|q|\nu(1+C_0).$ Then, we get that $\tilde{g}$ also satisfy the assumption (A2). Then by (A1), (A2) and (9), $(\tilde{Y}_s^{t+\varepsilon\wedge\tau},\tilde{Z}_s^{t+\varepsilon\wedge\tau})$ is a solution of BSDE with parameter $(\tilde{g},t+\varepsilon\wedge\tau,0)$ in $[t,t+\varepsilon\wedge\tau].$ By (3) and Lemma 3.1, there exists a constant $\gamma>0$ such that for each $\varepsilon\leq \gamma,$ we have
$$\sup_{t\leq s\leq t+\varepsilon\wedge\tau}|\tilde{Y}_{s}^{t+\varepsilon\wedge\tau}|\leq \tilde{b}\varepsilon e^{\alpha^+\varepsilon}\ \ \textrm{and} \ \ E\left[\int_{t}^{t+\varepsilon\wedge\tau}|\tilde{Z}_s^{t+\varepsilon\wedge\tau}|^2dr|{\cal{F}}_t\right]\leq \tilde{c}\varepsilon^2,\eqno(10)$$
where $\tilde{c}$ is a constant only depending on $\alpha, \beta, \tilde{b},T,\gamma.$\\

Set
\begin{eqnarray*}\
M^{\varepsilon,\tau}_t&:=&\frac{1}{\varepsilon}E\left[\int_{t}^{t+\varepsilon\wedge\tau}g(r,\tilde{Y}_r^{t+\varepsilon\wedge\tau}
+y+q\cdot(X_{r}^{t,x}-x),
\tilde{Z}_r^{t+\varepsilon\wedge\tau}+\sigma^\ast(r,X_{r}^{t,x})q)dr|{\cal{F}}_t\right]\\
P^{\varepsilon,\tau}_t&:=&\frac{1}{\varepsilon}E\left[\int_{t}^{t+\varepsilon\wedge\tau}g(r,y,
\sigma^\ast(r,x)q)dr|{\cal{F}}_t\right],\\
U^{\varepsilon,\tau}_t&:=&\frac{1}{\varepsilon}E\left[\int_{t}^{t+\varepsilon\wedge\tau}q\cdot b(r,X_{r}^{t,x})dr|{\cal{F}}_t\right],
\end{eqnarray*}
Then by (7) and (9), we have
\begin{eqnarray*}
&&\frac{1}{\varepsilon}\left(Y_t^{t+\varepsilon\wedge\tau}-y
\right)
-g(t,y,\sigma^\ast(t,x)q)-q\cdot b(t,x)\\
&=&\frac{1}{\varepsilon}\tilde{Y}_t^{t+\varepsilon\wedge\tau}
-g(t,y,\sigma^\ast(t,x)q)-q\cdot b(t,x)\\
&=&\left(M^{\varepsilon,\tau}_t-P^{\varepsilon,\tau}_t\right)
+\left(P^{\varepsilon,\tau}_t-g(t,y,\sigma^\ast(t,x)q)\right)+\left(U^{\varepsilon,\tau}_t-q\cdot b(t,x)\right).\ \ \ \  \ \ \ \  \ \  \ \ \ \  \ \  \ \ \ \  \ \  \ \ \ \ \ \ \ \ (11)
\end{eqnarray*}
By Jensen inequality, we have
\begin{eqnarray*}
 &&|M^{\varepsilon,\tau}_t-P^{\varepsilon,\tau}_t|\\
&\leq&\frac{1}{\varepsilon}E\left[\left|\int_{t}^{t+\varepsilon\wedge\tau}\left(g(r,\tilde{Y}_r^{t+\varepsilon\wedge\tau}+y
+q\cdot(X_{r}^{t,x}-x),
\tilde{Z}_r^{t+\varepsilon\wedge\tau}+\sigma^\ast(r,X_{r}^{t,x})q)-g(r,y,
\sigma^\ast(r,x)q)\right)dr\right||{\cal{F}}_t\right]\\
&\leq&\frac{1}{\varepsilon}E\left[\int_{t}^{t+\varepsilon\wedge\tau}\left|g(r,\tilde{Y}_r^{t+\varepsilon\wedge\tau}+y
+q\cdot(X_{r}^{t,x}-x),
\tilde{Z}_r^{t+\varepsilon\wedge\tau}+\sigma^\ast(r,X_{r}^{t,x})q)-g(r,y,
\sigma^\ast(r,x)q)\right|dr|{\cal{F}}_t\right]\\
\end{eqnarray*}
Then by (7), (8) and Lemma 3.2, we get there exists a non-negative process sequence $\{(\psi^n_t)_{t\in[0,T]}\}_{n=1}^{\infty}$ depending on $(y,x,q),$ which satisfies (5), such that for each $n\in \textbf{N}$
\begin{eqnarray*}
 &&|M^{\varepsilon,\tau}_t-P^{\varepsilon,\tau}_t|\\
&\leq&\frac{1}{\varepsilon}E\left[\int_{t}^{t+\varepsilon\wedge\tau}\left|2n\tilde{\lambda}_{2\tilde{M}}(|\tilde{Y}_r^{t+\varepsilon\wedge\tau}|
+|q||X_{r}^{t,x}-x|+|\tilde{Z}_r^{t+\varepsilon\wedge\tau}|^2+|X_{r}^{t,x}-x|^2)
+\psi_r^n\right|dr|{\cal{F}}_t\right]\\
&\leq&\frac{2}{\varepsilon}n\tilde{\lambda}_{2\tilde{M}}E\left[\int_{t}^{t+\varepsilon\wedge\tau}|\tilde{Y}_r^{t+\varepsilon\wedge\tau}|dr|{\cal{F}}_t\right]
+\frac{2}{\varepsilon}n\tilde{\lambda}_{2\tilde{M}}E\left[\int_{t}^{t+\varepsilon\wedge\tau}|\tilde{Z}_r^{t+\varepsilon\wedge\tau}|^2dr|{\cal{F}}_t\right]\\
&&+\frac{2}{\varepsilon}n\tilde{\lambda}_{2\tilde{M}}E\left[\int_{t}^{t+\varepsilon\wedge\tau}(|q||X_{r}^{t,x}-x|+|X_{r}^{t,x}-x|^2dr)|{\cal{F}}_t\right]
+\frac{1}{\varepsilon}E\left[\int_{t}^{t+\varepsilon}|\psi_r^n|dr|{\cal{F}}_t\right]\ \ \ \ \ \ \
\ \ \  (12)
\end{eqnarray*}
where $\tilde{\lambda}_{2\tilde{M}}=4\lambda_{2\tilde{M}} (1+|q|^2|\nu|^2),$ $\lambda_{2\tilde{M}}=\max\{|\alpha|, |\beta|, |\tilde{b}|,l(2\tilde{M})\}$. By (10), we can deduce
$$\lim\limits_{\varepsilon\rightarrow0^+}\frac{1}{\varepsilon}
E\left[\int_{t}^{t+\varepsilon\wedge\tau}|\tilde{Y}_r^{t+\varepsilon\wedge\tau}|dr|{\cal{F}}_t\right]
\leq\lim\limits_{\varepsilon\rightarrow0^+}\tilde{b}\varepsilon e^{\alpha^+\varepsilon}=0,\eqno(13)$$
and
$$\lim\limits_{\varepsilon\rightarrow0^+}\frac{1}{\varepsilon}
E\left[\int_{t}^{t+\varepsilon\wedge\tau}|\tilde{Z}_r^{t+\varepsilon\wedge\tau}|^2dr|{\cal{F}}_t\right]=0.\eqno(14)$$
By (6), Lebesgue dominated convergence theorem and the continuity of $X_{r}^{t,x}$ in $r$, we have,
\begin{eqnarray*}
\ \ \ \ \ \ \ \ \ \ \ \ \ \ \ \ \ \ \ &&\lim\limits_{\varepsilon\rightarrow0^+}\frac{1}{\varepsilon}E\left[\int_{t}^{t+\varepsilon\wedge\tau}(|q||X_{r}^{t,x}-x|+|X_{r}^{t,x}-x|^2)
dr|{\cal{F}}_t\right]\\
&=&E\left[\lim\limits_{\varepsilon\rightarrow0^+}\frac{1}{\varepsilon}\int_{t}^{t+\varepsilon\wedge\tau}(|q||X_{r}^{t,x}-x|+|X_{r}^{t,x}-x|^2)
dr|{\cal{F}}_t\right]\\
&=&0.\ \ \ \ \ \ \ \ \ \ \ \ \ \ \ \ \ \ \ \ \ \ \ \ \ \ \ \ \ \ \ \ \ \ \ \ \ \ \ \ \ \ \ \ \ \ \ \ \ \ \ \ \ \ \ \ \ \ \ \ \ \ \ \ \ \ \ \ \ \ \ \ \ \ \ \ \ \ \ \ \ \ \ \ \ \ \ \ \ \ (15)
\end{eqnarray*}
Since $\{(\psi^n_t)_{t\in[0,T]}\}_{n=1}^{\infty}$ satisfies (5), then by (12)-(15) and Lemma 3.3, we get that for almost every $t\in[0,T[$,
$$ P-a.s.,\ \ \lim\limits_{\varepsilon\rightarrow0^+}|M^{\varepsilon,\tau}_t
-P^{\varepsilon,\tau}_t|\leq\liminf_{n\rightarrow\infty}\lim\limits_{\varepsilon\rightarrow0^+}\frac{1}{\varepsilon}
E\left[\int_{t}^{t+\varepsilon}|\psi_r^n|dr|{\cal{F}}_t\right]
=\lim_{n\rightarrow\infty}|\psi_t^n|
=0.\eqno(16)$$
by (A2), (H2) and Lemma 3.3, we have for almost every $t\in[0,T[$,
$$P-a.s.,\ \ \lim\limits_{\varepsilon\rightarrow0^+}|P^{\varepsilon,\tau}_t-g(t,y,\sigma^\ast(t,x)q)|=0.\eqno(17)$$
By Jensen inequality, we have,
\begin{eqnarray*}
\ \ \ \ \ \ \  \ &&\lim\limits_{\varepsilon\rightarrow0^+}|U^{\varepsilon,\tau}_t-q\cdot b(t,x)|\\
&=&\lim\limits_{\varepsilon\rightarrow0^+}\left|E\left[\frac{1}{\varepsilon}\int_{t}^{t+\varepsilon\wedge\tau}
(q\cdot b(r,X_{r}^{t,x})-q\cdot b(r,x)+q\cdot b(r,x))dr|{\cal{F}}_t\right]-q\cdot b(t,x)\right|\\
&\leq&\lim\limits_{\varepsilon\rightarrow0^+}\left|\frac{1}{\varepsilon}E\left[\int_{t}^{t+\varepsilon\wedge\tau}(q\cdot b(r,X_{r}^{t,x})-q\cdot b(r,x))dr|{\cal{F}}_t\right]\right|\\
&&+\lim\limits_{\varepsilon\rightarrow0^+}\left|E\left[\frac{1}{\varepsilon}\int_{t}^{t+\varepsilon\wedge\tau}q\cdot b(r,x)dr|{\cal{F}}_t\right]-q\cdot b(t,x)\right|\\
&\leq&\lim\limits_{\varepsilon\rightarrow0^+}\frac{1}{\varepsilon}E\left[\int_{t}^{t+\varepsilon\wedge\tau}|q\cdot b(r,X_{r}^{t,x})-q\cdot b(r,x)|dr|{\cal{F}}_t\right]\\
&&+\lim\limits_{\varepsilon\rightarrow0^+}\left|E\left[\frac{1}{\varepsilon}\int_{t}^{t+\varepsilon\wedge\tau}q\cdot b(r,x)dr|{\cal{F}}_t\right]-q\cdot b(t,x)\right|. \ \ \ \ \ \ \ \ \ \ \ \ \ \ \ \ \ \ \ \ \ \ \ \ \ \ \ \ \ \ \ \ \ \ \ \ \ \ \ (18)
\end{eqnarray*}
By (H2) and Lemma 3.3, we have, for almost every $t\in[0,T[$,
$$P-a.s.,\ \ \lim\limits_{\varepsilon\rightarrow0^+}E\left[\frac{1}{\varepsilon}\int_{t}^{t+\varepsilon\wedge\tau}q\cdot b(r,x)dr|{\cal{F}}_t\right]-q\cdot b(t,x)=0.\eqno(19)$$
Then by (18), (19), (6), Lebesgue dominated convergence theorem, (H1) and the continuity of $X_{r}^{t,x}$ in $r$, we have, for almost every $t\in[0,T[$, $P-a.s.,$
\begin{eqnarray*}
\ \ \ \ \ \ \  \ \lim\limits_{\varepsilon\rightarrow0^+}|U^{\varepsilon,\tau}_t-q\cdot b(t,x)|&\leq&\lim\limits_{\varepsilon\rightarrow0^+}\frac{1}{\varepsilon}E\left[\int_{t}^{t+\varepsilon\wedge\tau}|q\cdot b(r,X_{r}^{t,x})-q\cdot b(r,x)|dr|{\cal{F}}_t\right]\\
&\leq&E\left[\lim\limits_{\varepsilon\rightarrow0^+}\frac{1}{\varepsilon}\int_{t}^{t+\varepsilon\wedge\tau}|q\cdot b(r,X_{r}^{t,x})-q\cdot b(r,x)|dr|{\cal{F}}_t\right]\\
&\leq&E\left[\lim\limits_{\varepsilon\rightarrow0^+}\frac{1}{\varepsilon}\int_{t}^{t+\varepsilon\wedge\tau}|q|
\mu|X_{r}^{t,x}-x|dr|{\cal{F}}_t\right]\\
&=&0.\ \ \ \ \ \ \ \ \ \ \  \ \ \ \ \ \ \ \ \ \ \ \ \ \ \ \ \ \ \ \ \ \ \ \ \ \ \ \ \ \ \ \ \ \ \ \ \ \ \ \ \ \ \ \ \ \ \ \ \ \ \ \ \ \ \ \ \ \ \ \ (20)
\end{eqnarray*}
By (11), (16), (17) and (20), The proof is complete.
$\Box$
\\

Now, by Theorem 3.4, we will give two representation theorems in $L^p$ and ${\cal{H}}^p$ for $p>0$, respectively.\\\\
\
\textbf{Corollary 3.5} Let $p>0$ and $g$ satisfy (A1) and (A2).  Then for each $(y,x,q)\in{\mathbf{R}}\times{\mathbf{R}}^{\mathit{m}}\times {\mathbf{R}}^{\mathit{m}}$ and almost every $t\in[0,T[$, we have
$$g\left(t,y,\sigma^\ast(t,x)q\right)+q\cdot b(t,x)
=L^p-\lim\limits_{\varepsilon\rightarrow0^+}\frac{1}{\varepsilon}
\left(Y_t^{t+\varepsilon\wedge\tau}-y\right),$$
where $\tau=\inf\{s\geq0:|X_{t+s}^{t,x}|> C_0\}\wedge (T-t)$ for a constant $C_0>|x|,$ and $(Y_s^{t+\varepsilon\wedge\tau},Z_s^{t+\varepsilon\wedge\tau})$ is an arbitrary solution of BSDE with parameter $(g,t+\varepsilon\wedge\tau,y+q\cdot(X_{t+\varepsilon\wedge\tau}^{t,x}-x)).$\\\\
\
\emph{Proof.}  By (7) and (10), for $\varepsilon\in]0,T-t]$, we have
$$P-a.s.,\ \ \frac{1}{\varepsilon}\left|Y_t^{t+\varepsilon\wedge\tau}-y\right|=\frac{1}{\varepsilon}|\tilde{Y}_{t}^{t+\varepsilon\wedge\tau}|\leq \tilde{b} e^{\alpha^+\varepsilon}\leq\tilde{b}e^{\alpha^+T}.\eqno(21)$$
For each $p>0,$ by (21), (A2), (H2) and Lebesgue dominated convergence theorem, we have
\begin{eqnarray*}
&&\lim\limits_{\varepsilon\rightarrow0^+}E\left|\frac{1}{\varepsilon}
\left(Y_t^{t+\varepsilon\wedge\tau}-y\right)-g\left(t,y,\sigma^\ast(t,x)q\right)+q\cdot b(t,x)\right|^p\\
&=&E\left|\lim\limits_{\varepsilon\rightarrow0^+}\frac{1}{\varepsilon}
\left(Y_t^{t+\varepsilon\wedge\tau}-y\right)-g\left(t,y,\sigma^\ast(t,x)q\right)+q\cdot b(t,x)\right|^p.
\end{eqnarray*}
Then by Theorem 3.4, we complete the proof.  $\Box$\\

Let $m=d, q=z, b(t,x)=0, \sigma(t,x)=1, x=0$ in Theorem 3.4. Then we have\\\\
\
\textbf{Corollary 3.6} Let $g$ satisfy the assumptions (A1) and (A2),  then for each $(y,z)\in{\mathbf{R}}\times{\mathbf{R}}^{\mathit{d}}$ and almost every $t\in[0,T[$, we have
$$P-a.s.,\ \ \ g\left(t,y,z\right)
=\lim\limits_{\varepsilon\rightarrow0^+}\frac{1}{\varepsilon}
\left(Y_t^{t+\varepsilon\wedge\tau}-y\right),$$
where $\tau=\inf\{s\geq0:|B_{t+s}-B_t|> C_0\}\wedge (T-t)$ for a constant $C_0>0$ and $(Y_s^{t+\varepsilon\wedge\tau},Z_s^{t+\varepsilon\wedge\tau})$ is an arbitrary solution of BSDE with parameter $(g,t+\varepsilon\wedge\tau,y+z\cdot(B_{t+\varepsilon\wedge\tau}-B_t)).$\\

In the end of this section, we will give the following remarks on the representation theorems obtained in the above.\\\\
\
\textbf{Remark 3.7}
\begin{itemize}
\item To our knowledge, all the representation theorems for BSDEs with linear growth are established in $L^p$ or ${\cal{H}}^p$ space for $1\leq p<2\ (\textrm{or}\ p=2)$ (see Briand et al. (2000), Jia (2008), Jiang (2005c, 2008), Fan and Jiang (2010), Fan et al. (2011), etc). But our representation theorems for quadratic BSDEs are established in $L^p$ sense for $p>0,$ due to the boundedness of solutions of such BSDEs.
\item To our knowledge, all the representation theorem for generators of BSDEs in $L^p$ space are established under the additional condition that $b(t,x)$ and $\sigma(t,x)$ in SDEs are right continuous in $t$ or independent on $t$ (see Briand et al. (2000), Jiang (2005), Jia (2008) and Jia and Zhang (2015), etc). But this condition is eliminated in our results.
\item Ma and Yao (2010, Theorem 4.1) established a representation theorem for quadratic BSDE using a different method under some assumptions slightly stronger than (A1)-(A4) and two additional assumptions (see conditions ($g1$) and ($g2$) in Ma and Yao (2010, Theorem 4.1)). Our representation theorems are not dependent on (A3), (A4) and such two additional assumptions.
\item Jia and Zhang (2015, Theorem 3.1) established a representation
theorem for BSDE with so called Lipschitz-quadratic generator (see assumption (H) in Jia and Zhang (2015)). In fact, one can check such Lipschitz-quadratic condition is a special case of (A1) and (A2).
\item The representation
theorems in Ma and Yao (2010, Theorem 4.1) and Jia and Zhang (2015, Theorem 3.1) both holds true for all $t\in[0,T[,$ due to the continuity assumption on generator $g.$ Without such continuity assumption, our representation theorem are obtained for almost every $t\in[0,T[$ due to the use of Lebesgue Lemma (Lemma 2.2).
\item For simplicity, we only consider the representation
theorem for quadratic BSDE with deterministic terminal time. In fact, with the similar method and some special treatments, one also can study the representation theorem for quadratic BSDE with random terminal time $\delta<\infty$, considered in Kobylanski (2000). But in this case, the representation theorem will hold true for almost every $t\in[0,\delta[.$
\end{itemize}
\section{Some applications for quadratic BSDEs}
In this section, we will give some applications of representation theorem for quadratic BSDEs.\\\\
{\textbf{Theorem 4.1}} (Converse comparison theorem) Let $g_1$ and $g_2$ satisfy (A1) and (A2), for any stopping time $\tau \in]0,T]$ and $\forall\xi\in L^\infty({\mathcal{F}}_\tau),$ BSDEs with parameter $(g_1,\tau,\xi)$ and $(g_2,\tau,\xi)$ exist solutions $(Y^{\tau,1},Z^{\tau,1})$ and $(Y^{\tau,2},Z^{\tau,2})$, respectively, such that $\forall t\in[0,T],$
$$P-a.s.,\ \ Y_{t\wedge\tau}^{\tau,1}\geq Y_{t\wedge\tau}^{\tau,2}.\eqno(22)$$
Then for each $(y,z)\in{\mathbf{R}}\times{\mathbf{R}}^{\mathit{d}}$ and almost every $t\in[0,T[$, we have
$$P-a.s.,\  \ g_1(t,y,z)\geq g_2(t,y,z).$$
\textit{Proof.} By Corollary 3.6, for each $(y,z)\in{\mathbf{R}}\times{\mathbf{R}}^{\mathit{d}}$ and almost every $t\in[0,T[$, there
exists a stopping time $\tau>0$, such that
$$P-a.s.,\ \ g_i(t,y,z)=\lim\limits_{n\rightarrow+\infty}n
(Y_t^{t+\frac{1}{n}\wedge\tau,i}-y),\ \ i=1,2. \eqno(23)$$
where $(Y_t^{t+\frac{1}{n}\wedge\tau,i},Z_t^{t+\frac{1}{n}\wedge\tau,i})$ is an arbitrary solution of BSDE with parameter $(g_i,t+\frac{1}{n}\wedge\tau,y+z\cdot(B_{t+\frac{1}{n}\wedge\tau}-B_t)),\ i=1,2$, respectively.
By (22) and (23), we can complete this proof.\ \  $\Box$\\

Self-financing condition and Zero-interest condition are considered in Jia (2008), Fan and Jiang (2010) and Fan et al. (2011) for BSDEs with linear growth. By Corollary 3.6, we can get the following similar results for quadratic BSDEs.
\\\\
\
{\textbf{Theorem 4.2} (Self-financing condition) Let $g$ satisfy (A1) and (A2). Then the following two statements are equivalent:

(i) For almost every $t\in[0,T[$,
$$P-a.s.,\ \ g(t,0,0)=0;$$

(ii) There exists a solution $(Y_t,Z_t)$ of BSDE with parameter $(g,T,0)$ such that $\forall t\in[0,T],$
$$P-a.s.,\ \ Y_t=0.$$
{\textbf{Theorem 4.3} (Zero-interest condition)} Let $g$ satisfy (A1) and (A2). Then the following two statements are equivalent:

(i) For each $y\in \textbf{R}$ and for almost every $t\in[0,T[$,
$$P-a.s.,\ \ g(t,y,0)=0;$$

(ii) For each $y\in \textbf{R},$ there exists a solution $(Y_t,Z_t)$ of BSDE with parameter $(g,T,y)$ such that $\forall t\in[0,T],$
$$P-a.s.,\ \ Y_t=y.$$

\section{Some applications for quadratic $g$-expectation}
In this section, using representation theorem obtained in this paper, we will give some properties of quadratic $g$-expectation in general case. Firstly we will give the following condition.
\begin{itemize}
  \item (A5). $P$-$a.s.$, $\forall (t,y)\in[0,T]\times{\mathbf{R}}^{\mathit{d}},\ \ g(t,y,0)=0.$
\end{itemize}

By Kobylanski (2000, Theorem 2.3 and Theorem 2.6), if $g$ satisfies (A1)-(A4), for each $\xi\in L^\infty({\mathcal {F}}_T),$ BSDE with parameter $(g,T,\xi)$ has a unique solution $(Y_t^T,Z_t^T)\in{\mathcal{S}}^\infty_T
({\mathbf{R}})\times{\cal{H}}^{2}_T({\mathbf{R}}^d)$. Moreover, if $g$ also satisfies (A5), we set ${\cal{E}}_g^T(\xi):=Y_0^T,$ called quadratic $g$-expectation of $\xi,$ and ${{\cal{E}}_g^T[\xi|{\cal{F}}_t]}=Y_t^T$ for $t\in[0,T],$ called conditional quadratic $g$-expectation of $\xi.$ For each stopping time $\sigma\in[0,T],$ We denote ${\cal{E}}_{g}^T[\xi|{\cal{F}}_{t\wedge\sigma}]:=Y_{t\wedge\sigma}^T,\ t\in[0,T].$\\\\
\
\textbf{Theorem 5.1} (Uniqueness theorem) Let $g_1$ and $g_2$ satisfy (A1)-(A5). Then the following two statements are equivalent:

(i) For each $\xi\in L^\infty({\mathcal {F}}_T),$ we have ${\cal{E}}_{g_1}^T(\xi)={\cal{E}}_{g_2}^T(\xi);$

(ii) For each $(y,z)\in{\mathbf{R}}\times{\mathbf{R}}^{\mathit{d}},$ and almost every $t\in[0,T[$, we have
$$P-a.s.,\  \ g_1(t,y,z)=g_2(t,y,z).$$
\emph{{Proof.}} We sketch this proof. By comparison theorem for quadratic BSDEs (see Kobylanski (2000, Theorem 2.6) or Ma and Yao (2010, Theorem 3.2)), we get (ii) $\Longrightarrow $ (i). One check the following fact: if (A1)-(A5) holds for $g,$ then for each stopping time $\sigma\in[0,T],\ \xi\in L^\infty({\mathcal {F}}_\sigma),$ we have
$$P-a.s.,\ \ {\cal{E}}_{g}^T[\xi|{\cal{F}}_{s\wedge\sigma}]={\cal{E}}_{g}^\sigma[\xi|{\cal{F}}_{s\wedge\sigma}], \ s\in[0,T].\eqno(24)$$
By the strict comparison in Ma and Yao (2010, Theorem 3.2) and the proof of Jiang (2004, Proposition 3.1), we also have the following fact: Let (A1)-(A5) holds for $g.$ If for each stopping time $\sigma\in[0,T]$ and $\forall\xi\in L^\infty({\mathcal {F}}_\sigma),$ we have ${\cal{E}}_{g_1}^\sigma[\xi]={\cal{E}}_{g_1}^\sigma[\xi],$ then
$\forall\xi\in L^\infty({\mathcal {F}}_\sigma)$ and for each stopping time $\rho\in[0,\sigma],$ we have $$P-a.s.,\ \ {\cal{E}}_{g_1}^\sigma[\xi|{\cal{F}}_{\rho}]={\cal{E}}_{g_1}^\sigma[\xi|{\cal{F}}_{\rho}].\eqno(25)$$
By (24), (25), Corollary 3.6 and the proof of Theorem 4.1, we can prove (i) $\Longrightarrow$ (ii). \ $\Box$\\

The following is a general converse comparison theorem for quadratic $g$-expectation, in which (i) $\Longleftrightarrow$ (iii) generalizes Ma and Yao (2010, Theorem 4.2) and (i) $\Longleftrightarrow$ (ii) is new. \\\\
\
\textbf{Theorem 5.2} (Converse comparison theorem) Let $g_1$ and $g_2$ satisfy (A1)-(A5). Then the following three statements are equivalent:

(i) For each $\xi\in L^\infty({\mathcal {F}}_T)$ and each $t\in[0,T],$ we have
$$P-a.s.,\ \ {\cal{E}}_{g_1}^T[\xi|{\cal{F}}_t]\leq{\cal{E}}_{g_2}^T[\xi|{\cal{F}}_t];$$

(ii) For each $\xi\in L^\infty({\mathcal {F}}_T)$ and each $t\in[0,T],$ we have
$$E\left[{\cal{E}}_{g_1}^T[\xi|{\cal{F}}_t]\right]\leq E\left[{\cal{E}}_{g_2}^T[\xi|{\cal{F}}_t]\right];$$

(iii) For each $(y,z)\in{\mathbf{R}}\times{\mathbf{R}}^{\mathit{d}}$ and almost every $t\in[0,T[$, we have
$$P-a.s.,\  \ g_1(t,y,z)\leq g_2(t,y,z).$$
\emph{Proof.} By comparison theorem for quadratic BSDEs (see Kobylanski (2000, Theorem 2.6) or Ma and Yao (2010, Theorem 3.2)), we can get (iii) $\Longrightarrow$ (i). By Theorem 4.1 and (24), we can get (i) $\Longrightarrow$ (iii). (i) $\Longrightarrow$ (ii) is trivial. So, we only need show (ii) $\Longrightarrow$ (i). In fact, If (i) does not hold, then there exist $\xi\in L^\infty({\mathcal {F}}_T)$ and $s\in[0,T]$ such that for some constant $\delta>0,$
$$P\left({\cal{E}}_{g_1}^T[\xi|{\cal{F}}_s]\geq{\cal{E}}_{g_2}^T[\xi|{\cal{F}}_s]+\delta\right)>0.$$
We set $A:=\left\{{\cal{E}}_{g_1}^T[\xi|{\cal{F}}_s]\geq{\cal{E}}_{g_2}^T[\xi|{\cal{F}}_s]+\delta\right\}.$ Clearly $A\in{\cal{F}}_s$, then by the properties of quadratic $g$-expectation (see Ma and Yao (2010)), we have
$${\cal{E}}_{g_1}^T[1_A\xi|{\cal{F}}_s]=1_A{\cal{E}}_{g_1}^T[\xi|{\cal{F}}_s]\geq1_A{\cal{E}}_{g_2}^T[\xi|{\cal{F}}_s]+1_A\delta
={\cal{E}}_{g_2}^T[1_A\xi|{\cal{F}}_s]+1_A\delta.$$
Taking expectation on both sides, we have
$$E\left[{\cal{E}}_{g_1}^T[1_A\xi|{\cal{F}}_s]\right]>E\left[{\cal{E}}_{g_2}^T[1_A\xi|{\cal{F}}_s]\right],$$
which contradicts (ii). The proof is complete. \ $\Box$\\

By Corollary 3.6 and the argument of Jiang (2008, Theorem 3.1) or Ma and Yao (2010, Proposition 4.3), we can get the following translation invariance of $g$-expectation, which generalizes Ma and Yao (2010, Proposition 4.3).\\\\
\
\textbf{Theorem 5.3} (Translation invariance) Let $g$ satisfy (A1)-(A5). Then the following three statements are equivalent:

(i) $g$ is independent on $y$;

(ii) For each $\xi\in L^\infty({\mathcal {F}}_T)$ and constant $C$, we have
$${\cal{E}}_{g}^T(\xi+C)={\cal{E}}_{g}^T(\xi)+C;$$

(iii) For each $\xi\in L^\infty({\mathcal {F}}_T),\ t\in[0,T]$ and $\eta\in{\cal{F}}_t,$ we have
$$P-a.s.,\ \ {\cal{E}}_{g}^T[\xi+\eta|{\cal{F}}_t]={\cal{E}}_{g}^T[\xi|{\cal{F}}_t]+\eta.$$

\end{document}